\documentclass[titlepage,draft,12pt]{article} 
\usepackage{amsfonts,latexsym,pstricks,pst-plot} 
\usepackage[a4paper]{geometry}

\newcommand{\n}{{\mathbb{N}}}
\newcommand{\cep}{c_{\ep}}

\newcommand{\rep}{r_{\ep}}
\newcommand{\roep}{\rho_{\ep}}
\newcommand{\uep}{u_{\ep}}

\newcommand{\Fep}{F_{\ep}}
\newcommand{\Gep}{G_{\ep}}

\newcommand{\tetep}{\theta_{\ep}}

\newcommand{\mdep}{\mathcal{D}_{\ep}}
\newcommand{\meep}{\mathcal{E}_{\ep}}
\newcommand{\mdepk}{\mathcal{D}_{\ep,k}}
\newcommand{\meepk}{\mathcal{E}_{\ep,k}}

\newcommand{\mGepk}{\mathcal{G}_{\ep,k}}
\newcommand{\auq}[1]{|A^{1/2}#1|^{2}}

\newcommand{\dau}{D(A^{1/2})}
\newcommand{\au}[1]{|A^{1/2}#1|}
\newcommand{\da}{D(A)}
\newcommand{\ep}{\varepsilon}

\newcommand{\m}[1]{m(\au{#1}^{2})}

\newtheorem{thm}{Theorem}[section]
\newtheorem{thmbibl}{Theorem}

\newtheorem{lemma}[thm]{Lemma}

\newtheorem{rmk}[thm]{Remark}

\title{Hyperbolic--parabolic singular perturbation for nondegenerate
Kirchhoff equations with critical weak dissipation} 
\author{Marina Ghisi\vspace{1ex}\\ {\normalsize
Universit\`a degli Studi di Pisa} \\{\normalsize Dipartimento di
Matematica ``Leonida Tonelli''}\\
{\normalsize 
PISA (Italy)}\\  
{\normalsize e-mail: \texttt{ghisi@dm.unipi.it}}\and 
Massimo Gobbino\vspace{1ex}\\ {\normalsize Universit\`a degli Studi di Pisa} 
\\{\normalsize Dipartimento di Matematica Applicata ``Ulisse Dini''}\\ 
{\normalsize 
 PISA (Italy)}\\  
{\normalsize e-mail: \texttt{m.gobbino@dma.unipi.it}}}
\date{}

\begin{document}
\maketitle
\begin{abstract}
	We consider the hyperbolic-parabolic singular perturbation problem
	for a nondegenerate quasilinear equation of Kirchhoff type with
	weak dissipation.  This means that the dissipative term is
	multiplied by a coefficient $b(t)$ which tends to 0 as $t\to
	+\infty$.
	
	The case where $b(t)\sim(1+t)^{-p}$ with $p<1$ has recently been
	considered.  The result is that the hyperbolic problem has a
	unique global solution, and the difference between solutions of
	the hyperbolic problem and the corresponding solutions of the
	parabolic problem converges to zero both as $t\to+\infty$ and as
	$\ep\to 0^{+}$.
	
	In this paper we show that these results cannot be true for $p>
	1$, but they remain true in the critical case $p=1$.
	
\vspace{1cm}

\noindent{\bf Mathematics Subject Classification 2000 (MSC2000):}
35B25, 35B40, 35L70.

\vspace{1cm} 
\end{abstract}
 
\section{Introduction}
Let $H$ be a real Hilbert space.  For every $x$ and $y$ in $H$, $|x|$
denotes the norm of $x$, and $\langle x,y\rangle$ denotes the scalar
product of $x$ and $y$.  Let $A$ be a self-adjoint linear operator on
$H$ with dense domain $D(A)$.  We assume that $A$ is nonnegative,
namely $\langle Ax,x\rangle\geq 0$ for every $x\in D(A)$, so that for
every $\alpha\geq 0$ the power $A^{\alpha}x$ is defined provided that
$x$ lies in a suitable domain $D(A^{\alpha})$.

Let $m:[0,+\infty)\to (0,+\infty)$ be a function of class $C^{1}$
satisfying the nondegeneracy condition
\begin{equation}
	\mu_{1}:=\inf_{\sigma\geq 0}m(\sigma)>0.
	\label{hp:ndg}
\end{equation}

We consider the second order Cauchy problem
\begin{equation}
	\ep\uep''(t)+b(t)\uep'(t)+\m{\uep(t)}A\uep(t)=0 \hspace{2em}
	\forall t\geq 0,
	\label{pbm:h-eq}
\end{equation}
\begin{equation}
	\uep(0)=u_0,\hspace{2em}\uep'(0)=u_1,
	\label{pbm:h-data}
\end{equation}
where $\ep>0$ is a parameter and $b:[0,+\infty)\to(0,+\infty)$ is
a given function.  We also consider the first order reduced Cauchy
problem
\begin{equation}
	b(t)u'(t)+\m{u(t)}Au(t)=0 \hspace{2em}\forall t\geq 0,
	\label{pbm:p-eq}
\end{equation}
\begin{equation}
	u(0)=u_0,
	\label{pbm:p-data}
\end{equation}
obtained setting formally $\ep=0$ in (\ref{pbm:h-eq}), and omitting
the second initial condition in~(\ref{pbm:h-data}).

It is well known that (\ref{pbm:h-eq}), (\ref{pbm:h-data}) is the
abstract setting of a nonlocal quasilinear hyperbolic partial
differential equation which was proposed as a model for the damped
small vibrations of an elastic string or membrane with uniform density
equal to $\ep$.  In the concrete setting assumption (\ref{hp:ndg}) is
equivalent to strict hyperbolicity.

Following the approach introduced by \textsc{J.\ L.\
Lions}~\cite{lions} in the linear case, we consider also the corrector
$\tetep(t)$ as the solution of the second order linear problem
\begin{equation}
	\ep\tetep''(t)+b(t)\tetep'(t)=0 \hspace{2em}
	\forall t\geq 0,
	\label{pbm:theta-eq}  
\end{equation}
\begin{equation}
	\tetep(0)=0,\hspace{2em}\tetep'(0)=u_1+\frac{1}{b(0)}
	\m{u_{0}}Au_{0}=:w_{0}.
	\label{pbm:theta-data}
\end{equation}

Since $\tetep'(0)=\uep'(0)-u'(0)$, this corrector keeps into account
the boundary layer due to the loss of one initial condition.  Finally
we define $\rep(t)$ and $\roep(t)$ in such a way that
\begin{equation}
	\uep(t)=u(t)+\tetep(t)+\rep(t)=u(t)+\roep(t)\quad\quad\forall
	t\geq 0.
	\label{defn:rep}
\end{equation}

The singular perturbation problem consists in proving that $\rep(t)\to
0$ or $\roep(t)\to 0$ in some sense as $\ep\to 0^{+}$.

\subparagraph{\emph{\textmd{Constant dissipation}}}

In the case of Kirchhoff equations this problem is well studied when
$b(t)$ is a positive constant.  The classical result (see
\cite{debrito}, \cite{yamada}) is the existence of a unique global
solution provided that $(u_{0},u_{1})\in\da\times\dau$ and $\ep$ is
small enough.  Existence of a global solution with $\ep$ large is
still an open problem, as well as the nondissipative case $b(t)\equiv
0$.

The behavior of solutions as $t\to +\infty$ and $\ep\to 0^{+}$ has
long been studied (see \cite{ew}, \cite{k-cattaneo},
\cite{gg:k-decay}, \cite{gg:k-PS}).  A complete answer was found by
\textsc{H.\ Hashimoto} and \textsc{T.\ Yamazaki} \cite{yamazaki}.
They proved that for initial data $(u_{0},u_{1})\in D(A^{3/2})\times
D(A)$ one has that
\begin{equation}
	|\roep(t)|+(1+t)^{1/2}|A^{1/2}\roep(t)|+
	\sqrt{\ep}(1+t)|\rep'(t)|\leq C\ep
	\quad\quad\forall t\geq 0,
	\label{history:c-diss}  
	\end{equation}
where of course $C$ doesn't depend on $t$ and $\ep$.  When
$(u_{0},u_{1})\in D(A^{2})\times D(A)$ the coefficient $\sqrt{\ep}$
may be dropped, thus providing a better convergence rate on
$\rep'(t)$.  We define (\ref{history:c-diss}) an \emph{error-decay
estimate} because it keeps into account in the same time both the
decay of solutions as $t\to +\infty$, and their behavior as $\ep\to
0^{+}$.  It unifies and extends or improves all previous convergence
or decay results.

This theory can be easily extended from the model case where $b(t)$ is
constant to the more general case where $b(t)$ is bounded by two positive
constants, with some extra conditions on $b'(t)$.

\subparagraph{\textmd{\emph{Subcritical weak dissipation}}}

The case where $b(t)\to 0$ as $t\to +\infty$, the typical example
being $b(t)=(1+t)^{-p}$ for some $p>0$, seems to be more difficult
because of the competition between the smallness of $\ep$ and the
smallness of $b(t)$.  To our knowledge the first results in this
direction were obtained by \textsc{K.\ Ono}~\cite{ono-wd}, who proved
global existence in the special degenerate case $m(\sigma)=\sigma$
provided that $p\leq 1/3$, and by \textsc{M.\ Nakao} and \textsc{J.\
Bae}~\cite{nakao} who proved global existence in the nondegenerate
case with a coercive operator and a special nonlinear dissipation term
which decays as $(1+t)^{-p}$ with $p<1$ as $t\to +\infty$.

Recently, \textsc{T.\ Yamazaki} \cite{yamazaki-wd} considered problem
(\ref{pbm:h-eq}), (\ref{pbm:h-data}) with the nondegeneracy assumption
(\ref{hp:ndg}) and $b(t)$ which decays as $(1+t)^{-p}$ with $p<1$ as
$t\to +\infty$.  She proved existence of a global solution when $\ep$
is small enough, and error-decay estimates for the singular
perturbation problem.  

The main result as stated in~\cite{yamazaki-wd} depends on several
parameters.  Limiting for simplicity to estimates such as
(\ref{history:c-diss}), she proved that for initial data
$(u_{0},u_{1})\in D(A^{5/2})\times D(A^{2})$ one has that
\begin{equation}
	|\roep(t)|+(1+t)^{(p+1)/2}|A^{1/2}\roep(t)|+
	\sqrt{\ep}(1+t)|\rep'(t)|\leq C\ep
	\quad\quad\forall t\geq 0,
	\label{history:w-diss}  
	\end{equation}
where of course $C$ doesn't depend on $t$ and $\ep$.  The coefficient
$\sqrt{\ep}$ may be dropped for more regular data.  This error-decay
estimate extends (\ref{history:c-diss}) (which becomes the special
case $p=0$) but for the fact that it seems to require more regularity
on the initial data.

In any case the time-decay rates on $|A^{1/2}\roep(t)|$ and
$|\rep'(t)|$ are those expected for solutions of the parabolic problem
(\ref{pbm:p-eq}), (\ref{pbm:p-data}), and they are optimal in the case
of noncoercive operators in the sense explained in \cite{ch}.  Roughly
speaking this means that for $p<1$ the smallness of $\ep$ is dominant
over the smallness of $b(t)$ as $t\to +\infty$, and in this regime
(\ref{pbm:h-eq}) behaves like a parabolic equation.

Proofs are based on a subtle spectral analysis of the corresponding
linearized equations which gives at the same time the existence result
and the decay-error estimates.

\subparagraph{\textmd{\emph{Supercritical and critical weak dissipation}}}

In this paper we consider problem (\ref{pbm:h-eq}), (\ref{pbm:h-data})
in the case where $b(t)=(1+t)^{-p}$ with $p\geq 1$.

When $p>1$ a simple argument (see Theorem~\ref{thm:p>1}) shows that
solutions of the hyperbolic problem (provided they globally exist,
which remains an open problem) do \emph{not} decay to 0 as $t\to
+\infty$.  On the contrary it is easy to see that the solutions of the
parabolic problem always decay to 0 as $t\to +\infty$ (faster and
faster as $p$ grows).  As a consequence one cannot expect
\emph{global-in-time} error-decay estimates such as
(\ref{history:w-diss}).  This also means that for $p>1$ the smallness
of $b(t)$ is dominant over the smallness of $\ep$, and in this regime
(\ref{pbm:h-eq}) behaves like a (nondissipative) hyperbolic equation.

This dichotomy, namely parabolic behavior for $p<1$ and hyperbolic
behavior for $p>1$, had already been observed in \cite{yamazaki-lin}
and \cite{wirth} in the case of linear equations (when $m(\sigma)$ is 
a positive constant).

Then we concentrate on the \emph{critical} case $p=1$, showing that
also in this case the parabolic nature prevails.  We prove indeed that
for every $(u_{0},u_{1})\in\da\times\dau$ the problem has a unique
global solution provided that $\ep$ is small enough, and this solution
decays to 0 as $t\to+\infty$ with the same rate of the corresponding
parabolic problem.  We also prove decay-error estimates for the
singular perturbation problem which extend (\ref{history:w-diss}) to
the case $p=1$.

Our approach is based uniquely on energy estimates and it applies
directly to the nonlinear problem.  The main advantage is of course
the possibility to treat the critical case.  Nevertheless we point out
that our assumptions on initial data are always minimal.  We obtain
indeed the global existence result for initial data in
$D(A)\times\dau$, which is of course the largest space where a
solution of class $C^{2}$ can be expected, and we obtain the
decay-error estimates for initial data in $D(A^{3/2})\times\dau$,
which is the largest space where error estimates of order $\ep$ can be
obtained, even when $b(t)$ and $m(\sigma)$ are positive constants
(see~\cite{gg:l-cattaneo}).  We also obtain the decay-error estimate
without the coefficient $\sqrt{\ep}$ for initial data
$(u_{0},u_{1})\in D(A^{2})\times D(A)$ (see Remark~\ref{rmk:nosqrt}).

For the sake of simplicity we work out the details only in the case
$p=1$.  On the other hand, the same technique applies to the case
$0<p<1$, or when $b(t)$ is constant, providing a different proof of
(\ref{history:c-diss}) and (\ref{history:w-diss}) with minimal
assumptions on the initial data.

When this paper was almost complete we were informed that Taeko
Yamazaki independently obtained some results on the same problem (see
\cite{yamazaki-cwd}) by different methods.

\setcounter{equation}{0}
\section{Statements}\label{sec:statements}

Throughout this paper the operator $A$ always satisfies the
following assumption:
\begin{list}{}{\leftmargin 4em \labelwidth 4em}
	\item[(Hp-$A$)] $A$ is a selfadjoint nonnegative linear operator
	whose domain $D(A)$ is dense in a Hilbert space $H$.  
\end{list}

We never assume $A$ to be coercive.

Existence of a unique global solution to problem (\ref{pbm:p-eq}),
(\ref{pbm:p-data}) is well known. The following is the result we need 
in this paper.

\begin{thmbibl}\label{thm:p-existence}
	Let $A$ be an operator satisfying (Hp-$A$), let $m:[0,+\infty)\to
	(0,+\infty)$ be a function of class $C^{1}$, let
	$b(t)=(1+t)^{-1}$, and let $u_{0}\in\da$.
	
	Then problem (\ref{pbm:p-eq}), (\ref{pbm:p-data}) has a unique
	solution
	$$u\in C^{1}([0,+\infty);H)\cap C^{0}([0,+\infty);\da).$$
	
	Moreover 
	$u\in C^{2}((0,+\infty);D(A^\alpha))$ for every $\alpha\geq 0$.
	
	If in addition $m$ satisfies (\ref{hp:ndg}), and $u_{0}\in
	D(A^{k/2})$ for some positive integer $k$, then there exists a
	constant $\gamma_{k}$, depending only on $k$ and $\mu_{1}$, such
	that
	\begin{equation}
		(1+t)^{2k}|A^{k/2}u(t)|^{2}\leq\gamma_{k}\left(
		|u_{0}|^{2}+|A^{k/2}u_{0}|^{2}\right)
		\quad\quad\forall t\geq 0,
		\label{th:parab-local}
	\end{equation}
	\begin{equation}
		\int_{0}^{+\infty}(1+s)^{2k+1}|A^{(k+1)/2}u(s)|^{2}\,ds\leq
		\gamma_{k}\left( |u_{0}|^{2}+|A^{k/2}u_{0}|^{2}\right).
		\label{th:parab-int}
	\end{equation}
\end{thmbibl}

The existence part of Theorem~\ref{thm:p-existence} is an easy
consequence of Theorem~4.1 in~\cite{k-par}.  In section~\ref{sec:par}
below we sketch the proof of the decay estimates
(\ref{th:parab-local}) and (\ref{th:parab-int}).

The first result of this paper concerns the global solvability of
problem (\ref{pbm:h-eq}), (\ref{pbm:h-data}), and the time decay of its
solutions.

\begin{thm}\label{thm:h-existence} 
	Let $A$ be an operator satisfying (Hp-$A$), let $m:[0,+\infty)\to
	(0,+\infty)$ be a function of class $C^{1}$ satisfying the
	non-degeneracy condition (\ref{hp:ndg}), let $b(t)=(1+t)^{-1}$,
	and let $(u_{0},u_{1})\in\da\times\dau$.
	
	Then there exists $\ep_{0}>0$ such that for every
	$\ep\in(0,\ep_{0})$ problem (\ref{pbm:h-eq}), (\ref{pbm:h-data})
	has a unique global solution 
	$$\uep\in C^{2}([0,+\infty);H)\cap C^{1}([0,+\infty);\dau) \cap
	C^{0}([0,+\infty);\da).$$
	
	Moreover there exists a constant $C$, independent on $\ep$ and
	$t$, such that
	\begin{equation}
		|\uep(t)|^{2}\leq C
		\quad\quad\forall t\geq 0;
		\label{th:decay-u}
	\end{equation}
	\begin{equation}
		\auq{\uep(t)}\leq \frac{C}{(1+t)^{2}}
		\quad\quad\forall t\geq 0;
		\label{th:decay-a1/2u}
	\end{equation}
	\begin{equation}
		|\uep'(t)|^{2}\leq \frac{C}{(1+t)^{2}}
		\quad\quad\forall t\geq 0;
		\label{th:decay-u'}
	\end{equation}
	\begin{equation}
		\int_{0}^{+\infty}(1+s)\left(
			|\uep'(s)|^{2}+\auq{\uep(s)}\right)\,ds\leq C;
		\label{th:decay-a1/2u-int}
	\end{equation}
	\begin{equation}
		\ep|A^{1/2}\uep'(t)|^{2}+ |A\uep(t)|^{2}\leq \frac{C}{(1+t)^{4}}
		\quad\quad\forall t\geq 0;
		\label{th:decay-au}
	\end{equation}
	\begin{equation}
		\int_{0}^{+\infty}(1+s)^{3}\left(
			|A^{1/2}\uep'(s)|^{2}+|A\uep(s)|^{2}\right)\,ds\leq C.
		\label{th:decay-au-int}
	\end{equation}
\end{thm}

The second result of this paper are the following decay-error
estimates.
\begin{thm}\label{thm:decay-error}
	Let $A$, $m$, $b$ be as in Theorem~\ref{thm:h-existence}.  Let
	$(u_{0},u_{1})\in D(A^{3/2})\times\dau$.  Let $\uep(t)$ be the
	solution of (\ref{pbm:h-eq}), (\ref{pbm:h-data}), let $u(t)$ be
	the solution of (\ref{pbm:p-eq}), (\ref{pbm:p-data}), and let
	$\rep(t)$ and $\roep(t)$ be defined by (\ref{defn:rep}).
	
	Then there exist $\ep_{0}>0$ and $C>0$ such that for every
	$\ep\in(0,\ep_{0})$ we have that
	\begin{equation}
		|\roep(t)|\leq C\ep
		\quad\quad\forall t\geq 0;
		\label{th:ps-u}
	\end{equation}
	\begin{equation}
		|A^{1/2}\roep(t)|\leq C\frac{\ep}{1+t}
		\quad\quad\forall t\geq 0;
		\label{th:ps-au}
	\end{equation}
	\begin{equation}
		|\rep'(t)|\leq C\frac{\sqrt{\ep}}{1+t}
		\quad\quad\forall t\geq 0;
		\label{th:ps-u'}
	\end{equation}
	\begin{equation}
		\int_{0}^{+\infty}(1+s)
		\left|A^{1/2}\roep(s)\right|^{2}\,ds\leq
		C\ep^{2};
		\label{th:ps-uint}
	\end{equation}
	\begin{equation}
		\int_{0}^{+\infty}(1+s)|\rep'(s)|^{2}\,ds\leq
		C\ep^{2}.
		\label{th:ps-u'int}
	\end{equation}
\end{thm}

The last result of this paper is that in the supercritical case $p>1$
solutions cannot decay to 0 as $t\to+\infty$.

\begin{thm}\label{thm:p>1}
	Let $A$ and $m$ be as in Theorem~\ref{thm:h-existence}. Let
	$b:[0,+\infty)\to(0,+\infty)$ be a continuous function such that
	\begin{equation}
		\int_{0}^{+\infty}b(s)\,ds<+\infty.
		\label{hp:int-b}
	\end{equation}
	
	Let $(u_{0},u_{1})\in\da\times\dau$ be such
	$|u_{1}|^{2}+|A^{1/2}u_{0}|^{2}>0$. Let us assume that for some
	$\ep>0$ problem (\ref{pbm:h-eq}), (\ref{pbm:h-data}) has a global 
	solution $\uep$.
	
	Then
	\begin{equation}
		\liminf_{t\to +\infty}\left(
		|\uep'(t)|^{2}+\auq{\uep(t)}\right)>0.
		\label{th:p>1}
	\end{equation}
\end{thm}

For the sake of simplicity we decided to state and prove our results
in the model case $b(t)=(1+t)^{-1}$ with minimal assumptions on
initial data.  The same technique however can be applied with more
regular initial data or more general dissipation terms.  In this
way one obtains the statements we mention in the two remarks below.

\begin{rmk}\label{rmk:nosqrt}
	\begin{em}
		Let us assume that in Theorem~\ref{thm:h-existence} we have
		that $(u_{0},u_{1})\in D(A^{(k+1)/2})\times D(A^{k/2})$ for
		some integer $k\geq 2$. Then it turns out that
		$$(1+t)^{2k}|A^{(k-1)/2}\uep'(t)|^{2}+(1+t)^{2k+2}\left(
		\ep|A^{k/2}\uep'(t)|^{2}+|A^{(k+1)/2}\uep(t)|^{2}\right) \leq
		C_{k},$$
		$$\int_{0}^{+\infty}(1+s)^{2k+1}\left(
		|A^{k/2}\uep'(s)|^{2}+|A^{(k+1)/2}\uep(s)|^{2}\right)\,ds\leq
		C_{k},$$
		where the constant $C_{k}$ depends on $k$, but not on $\ep$
		and $t$.  The decay rates are the same obtained in
		(\ref{th:parab-local}) and (\ref{th:parab-int}) for the
		parabolic problem.
		
		The decay-error estimates (\ref{th:ps-u}) through
		(\ref{th:ps-u'int}) admit similar extensions.  In particular
		when $(u_{0},u_{1})\in D(A^{(k+2)/2})\times D(A^{k/2})$ for
		some integer $k\geq 2$ we have that
		\begin{equation}
			(1+t)^{k-1}|A^{(k-2)/2}\rep'(t)|+
			(1+t)^{k}\left(\sqrt{\ep}|A^{(k-1)/2}\rep'(t)|+
			|A^{k/2}\roep(t)|\right)\leq C_{k}\ep,
			\label{th:rmk-1}
		\end{equation}
		\begin{equation}
			\int_{0}^{+\infty}(1+s)^{2k-1}\left(
			|A^{k/2}\roep(s)|^{2}+|A^{(k-1)/2}\rep'(s)|^{2}\right)\,ds\leq
			C_{k}\ep^{2}.
			\label{th:rmk-2}
		\end{equation}
		
		We sketch a proof of (\ref{th:rmk-1}) and (\ref{th:rmk-2}) in 
		section~\ref{sec:rmk}.
	\end{em}
\end{rmk}
\begin{rmk}
	\begin{em}
		Let us assume that $b(t)=1/\varphi(t)$, where
		$\varphi:[0,+\infty)\to(0,+\infty)$ is any function such that
		$\varphi'(t)$ is bounded, and
		$$\mbox{either }\int_{0}^{+\infty}
		\frac{[\varphi'(s)]^{2}}{\varphi(s)}\,ds<+\infty
		\quad\mbox{ or }\int_{0}^{+\infty}
		|\varphi''(s)|\,ds<+\infty.$$
		
		If $A$, $m$, $u_{0}$, $u_{1}$ are chosen as in
		Theorem~\ref{thm:h-existence}, and
		$$\Phi(t):=1+\int_{0}^{t}\varphi(s)\,ds,$$
		then we have that
		$$|\uep(t)|^{2}+\Phi(t)\left(|\uep'(t)|^{2}+
		\auq{\uep(t)}\right)+\Phi^{2}(t)\left(\ep|A^{1/2}\uep'(t)|^{2}+
		|A\uep(t)|^{2}\right)\leq C,$$
		$$\int_{0}^{+\infty}\varphi(s)\left[
		|\uep'(s)|^{2}+\auq{\uep(s)}+ \Phi(s)\left(
		|A^{1/2}\uep'(s)|^{2}+|A\uep(s)|^{2}\right)\right]\,ds\leq C.$$
		
		These estimates generalize (\ref{th:decay-u}) through
		(\ref{th:decay-au-int}).  Once again the decay rates of
		$\uep(t)$ are exactly those expected for $u(t)$.  The
		decay-error estimates (\ref{th:ps-u}) through
		(\ref{th:ps-u'int}) admit as well analogous extensions in
		terms of $\varphi(t)$ and $\Phi(t)$.
		
		We don't prove these estimates explicitly because this
		generality only complicates proofs without introducing new
		ideas.
	\end{em}
\end{rmk}

\setcounter{equation}{0}
\section{Proofs}\label{sec:proofs}

\subsection{Comparison results for ODEs}\label{sec:ODE}

In this section we state two comparison lemmata needed to prove our
main results.  We omit the simple standard proofs.

\begin{lemma}\label{lemma:ode-sqrt}
	Let $y:[0,+\infty)\to[0,+\infty)$ be a function of class $C^{1}$,
	and let $\psi:[0,+\infty)\to[0,+\infty)$ be a continuous
	function.  Let us assume that there exist two constants $c_{1}>0$
	and $c_{2}>0$ such that
	$$y'(t)\leq\psi(t)\left( -c_{1}y(t)+c_{2}\sqrt{y(t)}\right)
	\quad\quad\forall t\geq 0.$$
	
	Then $y(t)\leq\max\left\{y(0),(c_{2}/c_{1})^{2}\right\}$ for every
	$t\geq 0$.
\end{lemma}

\begin{lemma}\label{lemma:ode-lin}
	For $i=1,2,3$ let $g_{i}:[0,+\infty)\to[0,+\infty)$ be a
	continuous function.  Let $y:[0,+\infty)\to[0,+\infty)$ be a
	function of class $C^{1}$ such that $y(0)=0$ and
	$$y'(t)\leq -g_{1}(t)+g_{2}(t)y(t)+g_{3}(t)
	\quad\quad\forall t\geq 0.$$
	
	For $i=1,2,3$ let us set
	\begin{equation}
		G_{i}(t):=\int_{0}^{t}g_{i}(s)\,ds
		\quad\quad\forall t\geq 0.
		\label{defn:Gi}
	\end{equation}
	
	Then we have that
	$$y(t)+G_{1}(t)\leq e^{G_{2}(t)}G_{3}(t)
	\quad\quad\forall t\geq 0.$$
\end{lemma}

\subsection{Decay estimates for the parabolic equation}\label{sec:par}

For every $k\in\n$ and every $t>0$ we have that
\begin{eqnarray*}
	\left[(1+t)^{2k}|A^{k/2}u(t)|^{2}\right]' & = &
	2k(1+t)^{2k-1}|A^{k/2}u(t)|^{2}+  \\
	 &  & - 2(1+t)^{2k+1}\m{u(t)}|A^{(k+1)/2}u(t)|^{2},
\end{eqnarray*}
hence
\begin{eqnarray*}
	\lefteqn{\hspace{-5em}
	(1+t)^{2k}|A^{k/2}u(t)|^{2}+
	2\int_{0}^{t}(1+s)^{2k+1}\m{u(s)}|A^{(k+1)/2}u(s)|^{2}\,ds\ \leq}
	 \\
	 & \leq &
	 |A^{k/2}u_{0}|^{2}+2k\int_{0}^{t}(1+s)^{2k-1}|A^{k/2}u(s)|^{2}\,ds.
\end{eqnarray*}

Due to (\ref{hp:ndg}) this inequality implies (\ref{th:parab-local})
and (\ref{th:parab-int}) in the case $k=0$.  For $k>0$ the same
conclusions follow with an easy induction (we remind that all
intermediate norms $|A^{\alpha}u_{0}|$, with $0\leq\alpha\leq k/2$,
are controlled by $|u_{0}|$ and $|A^{k/2}u_{0}|$).

\subsection{Global existence and decay}\label{sec:existence}

In this section we prove Theorem~\ref{thm:h-existence}.

\subparagraph{\textmd{\emph{Local maximal solutions}}}

Problem (\ref{pbm:h-eq}), (\ref{pbm:h-data}) admits a unique local-in-time
solution, and this solution can be continued to a solution defined in 
a maximal interval $[0,T)$, where either $T=+\infty$, or
\begin{equation}
	\limsup_{t\to T^{-}}\left(|A^{1/2}\uep'(t)|^{2}+
	|A\uep(t)|^{2}\right)=+\infty.
	\label{limsup}
\end{equation}

We omit the proof of these standard results.  The interested reader is
referred to \cite{gg:k-dissipative} (see also \cite{AG}).

\subparagraph{\textmd{\emph{Monotonicity of the Hamiltonian}}}

Let
$$M(\sigma):=\int_{0}^{\sigma}m(s)\,ds,$$
and let 
\begin{equation}
	H(t):=\ep|\uep'(t)|^{2}+M(|A^{1/2}\uep(t)|^{2})
	\label{defn:H}
\end{equation}
be the usual Hamiltonian. From (\ref{pbm:h-eq}) we have that
\begin{equation}
	H'(t)=-2\frac{|\uep'(t)|^{2}}{1+t}
	\quad\quad\forall t\in[0,T),
	\label{deriv-H}
\end{equation}
hence
$$H(t)+2\int_{0}^{t}\frac{|\uep'(s)|^{2}}{1+s}\,ds=H(0)
\quad\quad\forall t\in[0,T).$$

Since $M(\sigma)\geq\mu_{1}\sigma$ this implies that
\begin{equation}
	|A^{1/2}\uep(t)|^{2}\leq\frac{H(0)}{\mu_{1}}
	\quad\quad\forall t\in[0,T).
	\label{est:auq-uep}
\end{equation}

\subparagraph{\textmd{\emph{Definitions and preliminaries}}}

Let us set for simplicity
\begin{equation}
	\cep(t):=\m{\uep(t)}.
	\label{defn:cep}
\end{equation}

Due to (\ref{hp:ndg}) and (\ref{est:auq-uep}) we have that
\begin{equation}
	\mu_{1}\leq\cep(t)\leq\mu_{2},
	\label{est:cep}
\end{equation}
where
\begin{equation}
	\mu_{2}:=\max\left\{m(\sigma):0\leq\sigma\leq\mu_{1}^{-1}H(0)\right\}.
	\label{defn:mu2}
\end{equation}

We also set
\begin{equation}
	L:=\max\left\{|m'(\sigma)|:0\leq\sigma\leq\mu_{1}^{-1}H(0)\right\},
	\label{defn:L}
\end{equation}
so that
\begin{equation}
	\frac{|\cep'(t)|}{\cep(t)}\leq\frac{2L}{\mu_{1}}\left|\langle
	u'(t),Au(t)\rangle\right|
	\quad\quad\forall t\in[0,T).
	\label{est:cep'}
\end{equation}

Let us consider the following constants, depending only on $u_{0}$,
$u_{1}$ (more precisely, on their norms in the spaces up to
$D(A)\times\dau$), $\mu_{1}$, $\mu_{2}$, $L$:
\begin{equation}
	k_{1}:=\max\left\{2,\frac{1}{\mu_{1}}\right\}\left(
	\frac{\mu_{2}}{\mu_{1}}\left(|u_{1}|^{2}+|u_{0}|^{2}\right)+
	|u_{1}|^{2}+\mu_{2}\auq{u_{0}}\right),
	\label{defn:k1}
\end{equation}
\begin{equation}
	k_{2}:=\max\left\{8,\frac{1}{\mu_{1}}\right\}
	\left(k_{1}+|u_{1}|^{2}+|u_{0}|^{2}\right),
	\label{defn:k2}
\end{equation}
\begin{equation}
	k_{3}:=k_{2}+\frac{1}{2}\left(\auq{u_{0}}+\auq{u_{1}}\right),
	\label{defn:k3}
\end{equation}
\begin{equation}
	k_{4}:=\max\left\{1,2\mu_{2}\right\}\left(
	\frac{|A^{1/2}u_{1}|^{2}}{\mu_{1}}+
	|Au_{0}|^{2}+\frac{4}{\mu_{1}}k_{3}\right).
	\label{defn:k4}
\end{equation}

Finally, let $\ep_{0}$ be small enough in such a way that
\begin{equation}
	\ep_{0}\leq\min\left\{\frac{1}{8},
	\frac{\mu_{1}}{8\mu_{2}}\right\},
	\hspace{2em}
	\frac{2L|\langle u_{1},Au_{0}\rangle|}{\mu_{1}}\ep_{0}<
	\frac{1}{2},
	\hspace{2em}
	\sqrt{\ep_{0}}\leq\frac{\mu_{1}}{2L(k_{1}+k_{4})}.
	\label{defn:ep0}
\end{equation}

The core of this proof are some estimates on the following energies:
\begin{eqnarray}
	D_{\ep,0}(t) & := & \frac{1-\ep}{2}|\uep(t)|^{2}+
	\ep(1+t)\langle\uep'(t),\uep(t)\rangle,
	\label{defn:D0}  \\
	D_{\ep,1}(t)& := &\frac{1-3\ep}{2}(1+t)^{2}|A^{1/2}\uep(t)|^{2}+
	\ep(1+t)^{3}\langle\uep'(t),A\uep(t)\rangle,
	\label{defn:D1}  \\
	\Fep(t)& := &\ep\frac{|A^{1/2}\uep'(t)|^{2}}{\cep(t)}+|A\uep(t)|^{2}.  
	\label{defn:F} \\
	\Gep(t)& := &(1+t)^{2}|\uep'(t)|^{2}.
	\label{defn:G}
\end{eqnarray}

\subparagraph{\textmd{\emph{Integral estimate on $A^{1/2}\uep$}}}

The time derivative of (\ref{defn:D0}) is
$$D_{\ep,0}'(t)=-(1+t)\cep(t)|A^{1/2}\uep(t)|^{2}+\ep(1+t)|\uep'(t)|^{2}.$$

Integrating in $[0,t]$ we obtain that
\begin{equation}
	\int_{0}^{t}(1+s)\cep(s)|A^{1/2}\uep(s)|^{2}\,ds=
	\ep\int_{0}^{t}(1+s)|\uep'(s)|^{2}\,ds+D_{\ep,0}(0)-D_{\ep,0}(t).
	\label{est:D0-int}
\end{equation}

Now we have that
$$D_{\ep,0}(0)=\frac{1-\ep}{2}|u_{0}|^{2}+\ep\langle u_{0},u_{1}\rangle
\leq\frac{1}{2}|u_{0}|^{2}+\frac{\ep}{2}|u_{1}|^{2},$$
and, since $\ep\leq 1/2$, we have that
$$-D_{\ep,0}(t)\leq-\frac{1}{4}|\uep(t)|^{2}+
\ep(1+t)|\uep'(t)|\cdot|\uep(t)|\leq
-\frac{1}{8}|\uep(t)|^{2}+2\ep^{2}(1+t)^{2}|\uep'(t)|^{2}.$$

Replacing these estimates in (\ref{est:D0-int}) we obtain that
\begin{eqnarray}
	\lefteqn{\hspace{-2em}\frac{1}{8}|\uep(t)|^{2}+
	\int_{0}^{t}(1+s)\cep(s)|A^{1/2}\uep(s)|^{2}\,ds\ \leq} 
	\nonumber \\
	 & \leq & \ep\int_{0}^{t}(1+s)|\uep'(s)|^{2}\,ds+
	2\ep^{2}(1+t)^{2}|\uep'(t)|^{2}
	+\frac{1}{2}|u_{0}|^{2}+\frac{\ep}{2}|u_{1}|^{2}.
	\label{est:D0-2}
\end{eqnarray}
	
\subparagraph{\textmd{\emph{First set of decay estimates}}}

We prove that
\begin{equation}
	(1+t)^{2}\ep|\uep'(t)|^{2}+(1+t)^{2}\auq{\uep(t)}+
	\int_{0}^{t}(1+s)|\uep'(s)|^{2}\,ds\leq k_{1},
	\label{decay:E}
\end{equation}
\begin{equation}
	|\uep(t)|^{2}+\int_{0}^{t}(1+s)|A^{1/2}\uep(s)|^{2}\,ds\leq
	k_{2},
	\label{est:D0-final}
\end{equation}
for every $t\in[0,T)$, where $k_{1}$ and $k_{2}$ are the constants
defined by (\ref{defn:k1}) and (\ref{defn:k2}).

To this end from (\ref{deriv-H}) we easily deduce that
\begin{equation}
	\left[(1+t)^{2}H(t)\right]' =
	-2(1-\ep)(1+t)|\uep'(t)|^{2}+2(1+t)M(\auq{\uep(t)}).
	\label{eq:deriv-E}
\end{equation}

Moreover from (\ref{est:auq-uep}), (\ref{est:cep}), and (\ref{defn:mu2})
we have that
$$M(\auq{\uep(t)})\leq\mu_{2}\auq{\uep(t)}\leq
\frac{\mu_{2}}{\mu_{1}}\cep(t)\auq{\uep(t)}.$$ 

This allows to estimate the last term in (\ref{eq:deriv-E}), yielding
that $$\left[(1+t)^{2}H(t)\right]'\leq-2(1-\ep)(1+t)|\uep'(t)|^{2}+
\frac{2\mu_{2}}{\mu_{1}} (1+t)\cep(t)\auq{\uep(t)}.$$

Now we integrate in $[0,t]$ and we use  (\ref{est:D0-2}). After
rearranging the terms we obtain that
$$(1+t)^{2}\left(1-\frac{4\mu_{2}}{\mu_{1}}\ep\right)\ep|\uep'(t)|^{2}
+(1+t)^{2}M(|A^{1/2}\uep(t)|^{2})+$$
$$+2\left(1-\ep-\frac{\mu_{2}}{\mu_{1}}\ep\right)
\int_{0}^{t}(1+s)|\uep'(s)|^{2}\,ds\leq H(0)+
\frac{\mu_{2}}{\mu_{1}}\left(|u_{0}|^{2}+\ep|u_{1}|^{2}\right).$$

From our smallness assumptions on $\ep$ and the usual estimates on
$M(\auq{\uep(t)})$ we deduce that
$$\frac{1}{2}(1+t)^{2}\ep|\uep'(t)|^{2}
+\mu_{1}(1+t)^{2}|A^{1/2}\uep(t)|^{2}+
\int_{0}^{t}(1+s)|\uep'(s)|^{2}\,ds\ \leq$$
$$\leq\ |u_{1}|^{2}+\mu_{2}\auq{u_{0}}+
\frac{\mu_{2}}{\mu_{1}}\left(|u_{0}|^{2}+|u_{1}|^{2}\right),$$
from which (\ref{decay:E}) easily follows.  

This allows to estimate the first two terms in the right-hand side of
(\ref{est:D0-2}).  We thus obtain (\ref{est:D0-final}).

\subparagraph{\textmd{\emph{Integral estimate on $A\uep$}}}

The time derivative of (\ref{defn:D1}) is
$$D_{\ep,1}'(t)=-(1+t)^{3}\cep(t)|A\uep(t)|^{2}+
\ep(1+t)^{3}\auq{\uep'(t)}+
(1-3\ep)(1+t)\auq{\uep(t)}.$$

Integrating in $[0,t]$ we obtain that
\begin{eqnarray}
	\int_{0}^{t}(1+s)^{3}\cep(s)|A\uep(s)|^{2}\,ds & \leq & 
	\ep\int_{0}^{t}(1+s)^{3}\auq{\uep'(s)}\,ds+
	D_{\ep,1}(0)+
	\nonumber  \\
	 &  & +\int_{0}^{t}(1+s)\auq{\uep(s)}\,ds-D_{\ep,1}(t).
	\label{est:D1-int}
\end{eqnarray}

The second integral in the right-hand side can be estimated using
(\ref{est:D0-final}). Moreover
$$D_{\ep,1}(0)=\frac{1-3\ep}{2}|A^{1/2}u_{0}|^{2}+\ep\langle
A^{1/2}u_{0},A^{1/2}u_{1}\rangle
\leq\frac{1}{2}|A^{1/2}u_{0}|^{2}+\frac{\ep}{2}|A^{1/2}u_{1}|^{2},$$
and, since $\ep\leq 1/6$, we have that
\begin{eqnarray*}
	-D_{\ep,1}(t) &\leq  & -\frac{1}{4}(1+t)^{2}|A^{1/2}\uep(t)|^{2}+
	2\cdot(1+t)\frac{|A^{1/2}\uep(t)|}{2}\cdot
	\ep(1+t)^{2}|A^{1/2}\uep'(t)| \\
	 & \leq & \ep^{2}(1+t)^{4}|A^{1/2}\uep'(t)|^{2}.
\end{eqnarray*}

Replacing all these estimates in (\ref{est:D1-int}) we obtain that
\begin{eqnarray}
	\lefteqn{\hspace{-2em}
	\int_{0}^{t}(1+s)^{3}\cep(s)|A\uep(s)|^{2}\,ds\ \leq} 
	\nonumber \\
	 & \leq & \ep\int_{0}^{t}(1+s)^{3}\auq{\uep'(s)}\,ds+
	 \ep^{2}(1+t)^{4}\auq{\uep'(t)}+k_{3}.
	\label{est:D1-2}
\end{eqnarray}

\subparagraph{\textmd{\emph{Second set of decay estimates}}}

Let us set
\begin{equation}
	S:=\sup\left\{\tau\leq T:\ep\frac{|\cep'(t)|}{\cep(t)}\leq
	\frac{1}{2(1+t)}\quad\forall t\in[0,\tau]\right\}.
	\label{defn:S}
\end{equation}

We claim that $S>0$ and
\begin{equation}
	(1+t)^{4}\ep|A^{1/2}\uep'(t)|^{2}+(1+t)^{4}|A\uep(t)|^{2}+
	\int_{0}^{t}(1+s)^{3}|A^{1/2}\uep'(s)|^{2}\,ds\leq k_{4},
	\label{est:F}
\end{equation}
\begin{equation}
	\int_{0}^{t}(1+s)^{3}|A\uep(s)|^{2}\,ds\leq
	\frac{k_{4}+k_{3}}{\mu_{1}},
	\label{est:D1-final}
\end{equation}
for every $t\in[0,S)$, where $k_{3}$ and $k_{4}$ are the constants
defined by (\ref{defn:k3}) and (\ref{defn:k4}).

Let us prove these claims. Thanks to inequality (\ref{est:cep'}) with 
$t=0$ and the second inequality in (\ref{defn:ep0}) we have that
$$\ep\frac{|\cep'(0)|}{\cep(0)}\leq\ep_{0}
\frac{2L|\langle u_{1},Au_{0}\rangle|}{\mu_{1}}<\frac{1}{2},$$
hence $S>0$. Moreover we have that
$$\left[(1+t)^{4}\Fep(t)\right]' = -(1+t)^{4}
\left(\frac{2-4\ep}{1+t}+\ep\frac{\cep'(t)}{\cep(t)}\right)
\frac{\auq{\uep'(t)}}{\cep(t)}+4(1+t)^{3}|A\uep(t)|^{2}$$
for every $t\in[0,T)$.  Now let us restrict to $t\in[0,S)$.  Due to
the definition of $S$ and the fact that $\ep\leq 1/8$ this implies
that 
$$\left[(1+t)^{4}\Fep(t)\right]' \leq -(1+t)^{3}
\frac{\auq{\uep'(t)}}{\cep(t)}+\frac{4}{\mu_{1}}(1+t)^{3}\cep(t)
|A\uep(t)|^{2}.$$

Now we integrate in $[0,t]$ and we use (\ref{est:D1-2}). After
recollecting the terms we end up with 
\begin{eqnarray}
	\lefteqn{\hspace{-23em}
	\left(\frac{1}{\cep(t)}-\frac{4}{\mu_{1}}\ep\right)
	(1+t)^{4}\ep|A^{1/2}\uep'(t)|^{2} +(1+t)^{4}|A\uep(t)|^{2}+}
	\nonumber\\
	+\int_{0}^{t}\left(\frac{1}{\cep(s)}-\frac{4}{\mu_{1}}\ep\right)
	(1+s)^{3}|A^{1/2}\uep'(s)|^{2}\,ds & \leq & F(0)+
	\frac{4}{\mu_{1}}k_{3}.
	\label{est:F-interm}
\end{eqnarray}

Since $\ep\leq\mu_{1}/(8\mu_{2})$ we have that
$$\frac{1}{\cep(t)}-\frac{4}{\mu_{1}}\ep\geq
\frac{1}{\mu_{2}}-\frac{4}{\mu_{1}}\ep\geq
\frac{1}{2\mu_{2}},$$
and therefore (\ref{est:F}) easily follows from (\ref{est:F-interm}).
This in turn allows to estimate the first two terms in the right-hand
side of (\ref{est:D1-2}).  We thus obtain (\ref{est:D1-final}).

\subparagraph{\textmd{\emph{Global existence}}}

We prove that $S=T=+\infty$.  Let us assume by contradiction that
$S<T$.  By definition (\ref{defn:S}) of $S$ this means that
necessarily
\begin{equation}
	\ep\frac{|\cep'(S)|}{\cep(S)}=\frac{1}{2(1+S)}
	\label{eq:S-lim}
\end{equation}

On the other hand, applying (\ref{est:cep'}) with $t=S$,
(\ref{decay:E}), (\ref{est:F}), and the last inequality in
(\ref{defn:ep0}), we have that
$$\ep\frac{|\cep'(S)|}{\cep(S)}\; \leq\; 
\ep\frac{2L}{\mu_{1}}\left|\langle
\uep'(S),A\uep(S)\rangle\right|\;\leq\;
\frac{L}{\mu_{1}}\sqrt{\ep}\left(
\ep|\uep'(S)|^{2}+|A\uep(S)|^{2}\right)\; \leq$$
$$\leq\; \frac{L}{\mu_{1}}\sqrt{\ep_{0}}\left(
\frac{k_{1}}{(1+S)^{2}}+\frac{k_{4}}{(1+S)^{4}}\right)\; <\;
\frac{1}{2(1+S)},$$
which contradicts (\ref{eq:S-lim}).

It remains to prove that $T=+\infty$. To this end it is enough to show
that (\ref{limsup}) cannot be true, and this immediately follows from 
(\ref{est:F}).

\subparagraph{\textmd{\emph{Last decay estimate}}}

Since $T=+\infty$ we know that (\ref{decay:E}) and
(\ref{est:D0-final}) hold true for every $t\geq 0$.  This proves
(\ref{th:decay-u}), (\ref{th:decay-a1/2u}), and
(\ref{th:decay-a1/2u-int}).  Since $S=+\infty$ we know that also
(\ref{est:F}) and (\ref{est:D1-final}) hold true for every $t\geq
0$.  This proves (\ref{th:decay-au}) and (\ref{th:decay-au-int}).

It remains to prove (\ref{th:decay-u'}).  To this end we compute the
time derivative of (\ref{defn:G}):
$$\Gep'(t)=-2\left(\frac{1}{\ep}-1\right)(1+t)|\uep'(t)|^{2}-
\frac{2}{\ep}(1+t)^{2}\langle\uep'(t),\cep(t)A\uep(t)\rangle.$$

From estimate (\ref{est:F}) and the fact that $\ep\leq 1/2$ 
it follows that
\begin{eqnarray*}
	\Gep'(t) & \leq & -\frac{\Gep(t)}{\ep(1+t)}+\frac{2\mu_{2}}{\ep}
	(1+t)|A\uep(t)|\sqrt{\Gep(t)} \\
	 & \leq & \frac{1}{\ep(1+t)}\left(-\Gep(t)+2\mu_{2}\sqrt{k_{4}}
	\sqrt{\Gep(t)}\right).
\end{eqnarray*}

Therefore the conclusion follows from Lemma~\ref{lemma:ode-sqrt}
applied with $y(t):=\Gep(t)$.

\subsection{Decay-error estimates}\label{sec:ps}

In this section we prove Theorem~\ref{thm:decay-error}.

\subparagraph{\textmd{\emph{Notations and preliminaries}}}

Throughout this proof we assume that $\ep\in(0,\ep_{0})$, where
$\ep_{0}$ satisfies inequalities (\ref{defn:ep0}) as in the proof of the global
existence result, and the further inequality
\begin{equation}
	\ep_{0}\leq\frac{\mu_{1}}{128\mu_{2}}.
	\label{defn:ep0-bis}
\end{equation}

Under assumption (\ref{defn:ep0}) we already know that problem
(\ref{pbm:h-eq}), (\ref{pbm:h-data}) has a unique global solution
satisfying (\ref{th:decay-u}) through (\ref{th:decay-au-int}), and
\begin{equation}
	\ep\frac{|\cep'(t)|}{\cep(t)}\leq\frac{1}{2(1+t)}
	\quad\quad\forall t\geq 0,
	\label{est:cep'-cep}
\end{equation}
where of course $\cep(t)$ is defined by (\ref{defn:cep}).  Accordingly
we set $c(t):=\m{u(t)}$.  The function $t\to M(\auq{u(t)})$ turns out
to be nonincreasing, and therefore inequalities (\ref{est:auq-uep}) and
(\ref{est:cep}) hold true also with $u(t)$ and $c(t)$ in place of
$\uep(t)$ and $\cep(t)$.

The corrector $\tetep(t)$, namely the solution of
(\ref{pbm:theta-eq}), (\ref{pbm:theta-data}), can be explicitly
computed to be
\begin{equation}
	\tetep(t)=\frac{\ep}{1-\ep}\left(1-
	(1+t)^{1-1/\ep}\right)w_{0},
	\label{defn:tetep}
\end{equation} 
so that
\begin{equation}
	\tetep'(t)=w_{0}(1+t)^{-1/\ep}.
	\label{eq:tetep'}
\end{equation} 

Simple calculations show that $\rep$ is the
solution of the Cauchy problem
$$\ep\rep''(t)+\frac{1}{1+t}\rep'(t)+\cep(t)
A\roep(t)=(c(t)-\cep(t))Au(t)-\ep u''(t),$$
$$\rep(0)=0,\hspace{2em}\rep'(0)=0,$$
while $\roep$ is the solution of the Cauchy problem
$$\ep\roep''(t)+\frac{1}{1+t}\roep'(t)+\cep(t)
A\roep(t)=(c(t)-\cep(t))Au(t)-\ep u''(t),$$
$$\roep(0)=0,\hspace{2em}\roep'(0)=w_{0}.$$

In order to estimate $\roep$ and $\rep$ we introduce the energies
\begin{equation}
	\mdep(t):=\frac{1-\ep}{2}|\roep(t)|^{2}+\ep(1+t)
	\langle\roep'(t),\roep(t)\rangle,
	\label{defn:mdep}
\end{equation}
\begin{equation}
	\meep(t):=\ep\frac{|\rep'(t)|^{2}}{\cep(t)}+\auq{\roep}.
	\label{defn:meep}
\end{equation}

From now on $K_{1}$, $K_{2}$, $\alpha_{1}$, \ldots, $\alpha_{23}$
denote constants depending only on $u_{0}$, $u_{1}$ (more precisely on
their norms in the spaces up to $D(A^{3/2})\times D(A^{1/2})$),
$\mu_{1}$, $\mu_{2}$, $L$.

\subparagraph{\textmd{\emph{Integral estimate on  $u''$}}}

We prove that when $u_{0}\in D(A^{3/2})$ we have that
\begin{equation}
	\int_{0}^{+\infty}(1+s)^{3}|u''(s)|^{2}\,ds<+\infty.
	\label{est:u''-int}
\end{equation}

Indeed for every $t>0$ we have that
$$u''(t)=-c(t)Au(t)+(1+t)^{2}c^{2}(t)A^{2}u(t)+
2(1+t)^{2}c(t)m'(|A^{1/2}u(t)|^{2})|Au(t)|^{2}Au(t).$$

From (\ref{th:parab-local}) with $k=1$ we deduce that $|A^{1/2}u(t)|$
is bounded, hence also $c(t)$ and $m'(|A^{1/2}u(t)|^{2})$ are bounded.
Moreover also $ (1+t)^{2}|Au(t)|^{2}$ is bounded because of
(\ref{th:parab-local}) with $k=2$. It follows that
$$|u''(t)|\leq \alpha_{1}|Au(t)|+\alpha_{2}(1+t)^{2}|A^{2}u(t)|,$$
hence 
$$(1+t)^{3}|u''(t)|^{2}\leq \alpha_{3}(1+t)^{3}|Au(t)|^{2}+
\alpha_{4}(1+t)^{7}|A^{2}u(t)|^{2}.$$

Therefore (\ref{est:u''-int}) follows from (\ref{th:parab-int})
applied with $k=1$ and $k=3$.

\subparagraph{\textmd{\emph{Boundedness of $\meep$}}}

We prove that there exists a constant $K_{1}$ such that
\begin{equation}
	\meep(t)+\int_{0}^{t}\frac{1}{1+s}\frac{|\rep'(s)|^{2}}{\cep(s)}\,ds
	\leq K_{1}\ep^{2},
	\hspace{2em}
	\forall t\geq 0.
	\label{est:meep-1}
\end{equation}

The time derivative of (\ref{defn:meep}) is
\begin{eqnarray}
	\meep'(t) & = & -\frac{|\rep'(t)|^{2}}{\cep(t)}\left(
	\ep\frac{\cep'(t)}{\cep(t)}+\frac{2}{1+t}\right)+
	2\langle\tetep'(t),A\roep(t)\rangle+ 
	\nonumber  \\
	 &  & +\frac{2}{\cep(t)}\langle\rep'(t),(c(t)-\cep(t))Au(t)-\ep
	 u''(t)\rangle.
	 \label{eq:meep'}
\end{eqnarray}

Let us estimate the three terms in the right-hand side. The first
summand can be easily estimated by (\ref{est:cep'-cep}). For the
second summand we use (\ref{eq:tetep'}) and
we obtain that
$$2\langle\tetep'(t),A\roep(t)\rangle\leq
2|A^{1/2}\tetep'(t)|\cdot|A^{1/2}\roep(t)|\leq
2|A^{1/2}w_{0}|(1+t)^{-1/\ep}|A^{1/2}\roep(t)|.$$

From (\ref{defn:L}), (\ref{th:parab-local}) with $k=1$, and
(\ref{th:decay-a1/2u}) we have that
\begin{eqnarray}
	|c(t)-\cep(t)| & \leq & L\left|
	\auq{u(t)}-\auq{\uep(t)}\right|  
	\nonumber \\
	 & \leq & L\left(
	\au{u(t)}+\au{\uep(t)}\right)\au{\roep(t)} 
	\nonumber  \\
	 & \leq & \alpha_{5}\au{\roep(t)},
	 \label{est:cep-c}
\end{eqnarray}
hence (for shortness' sake we omit the dependence on $t$)
\begin{eqnarray*}
	 \frac{2}{\cep}\langle\rep',(c-\cep)Au-\ep
	 u''\rangle & \leq & \frac{1}{\cep}\cdot 2
	 \frac{|\rep'|}{\sqrt{1+t}}\cdot
	 \sqrt{1+t}\left(|c-\cep||Au|+\ep |u''|\right)  \\
	 & \leq & \frac{1}{2}\frac{|\rep'|^{2}}{(1+t)\cep}+
	 \frac{4}{\cep}(1+t)\left(|c-\cep|^{2}|Au|^{2}+
	 \ep^{2} |u''|^{2}\right)  \\
	 & \leq & \frac{1}{2}\frac{|\rep'|^{2}}{(1+t)\cep}+
	 \alpha_{6}(1+t)\auq{\roep}|Au|^{2}+\alpha_{7}
	 \ep^{2}(1+t) |u''|^{2}.
\end{eqnarray*}

Replacing all these estimates in (\ref{eq:meep'}) we obtain
that 
\begin{eqnarray}
	\meep'(t) & \leq & -\frac{1}{1+t}\frac{|\rep'(t)|^{2}}{\cep(t)}+
	\alpha_{6}(1+t)|Au(t)|^{2}\auq{\roep(t)}+ 
	\nonumber \\
	&&
	+\alpha_{7}(1+t)\ep^{2}|u''(t)|^{2}+
	\alpha_{8}(1+t)^{-1/\ep}\au{\roep(t)},
	\label{est:meep'}  
\end{eqnarray}
and in particular
\begin{eqnarray*}
	\meep'(t)  & \leq & -\frac{1}{1+t}\frac{|\rep'(t)|^{2}}{\cep(t)}+
	\alpha_{6}(1+t)|Au(t)|^{2}\meep(t)+  \\
	&&  +\alpha_{7}(1+t)\ep^{2}|u''(t)|^{2}+
	\alpha_{8}(1+t)^{-1/\ep}\sup_{t\geq 0}\au{\roep(t)}.  \\
	& =: & -g_{1}(t)+g_{2}(t) \meep(t)+g_{3}(t).
\end{eqnarray*}

We can now apply Lemma~\ref{lemma:ode-lin} with $y(t):=\meep(t)$ (we
recall that $\meep(0)=0$). The function $G_{2}(t)$ defined according
to (\ref{defn:Gi}) is bounded in $t$ because of (\ref{th:parab-int})
with $k=1$. It follows that
\begin{eqnarray*}
	\meep(t)+\int_{0}^{t}\frac{1}{1+s}\frac{|\rep'(s)|^{2}}{\cep(s)}\,ds
	& \leq & \alpha_{9}\ep^{2}\int_{0}^{t}(1+s)|u''(s)|^{2}\,ds+ 
	\\
	 & & +\alpha_{9}\sup_{t\geq
	 0}\au{\roep(t)}\int_{0}^{t}(1+s)^{-1/\ep}\,ds.
\end{eqnarray*}

The first integral in the right-hand side is uniformly bounded due to
(\ref{est:u''-int}), the second one is less than $2\ep$.  It follows that
\begin{equation}
	\meep(t)+\int_{0}^{t}\frac{1}{1+s}\frac{|\rep'(s)|^{2}}{\cep(s)}\,ds
	\leq \alpha_{10}\ep^{2}+\alpha_{11}\ep\sup_{t\geq 0}\au{\roep(t)}.
	\label{est:meep-0}
\end{equation}

In particular
$$\left(\sup_{t\geq 0}\au{\roep(t)}\right)^{2}\leq
\alpha_{10}\ep^{2}+\alpha_{11}\ep\sup_{t\geq 0}\au{\roep(t)}\leq
\frac{1}{2}\left(\sup_{t\geq 0}\au{\roep(t)}\right)^{2}+
\alpha_{12}\ep^{2},$$
and therefore
\begin{equation}
	\sup_{t\geq 0}\au{\roep(t)}\leq \alpha_{13}\ep.
	\label{est:sup}
\end{equation}

Coming back to (\ref{est:meep-0}) we obtain (\ref{est:meep-1}).

\subparagraph{\textmd{\emph{Integral estimate on $A^{1/2}\roep$}}}

The time derivative of (\ref{defn:mdep}) is
\begin{eqnarray*}
	\mdep'(t) & = &
	-(1+t)\cep(t)\auq{\roep(t)}+\ep(1+t)|\roep'(t)|^{2}+  \\
	 & & +(1+t)(c(t)-\cep(t))\langle Au(t),\roep(t)\rangle- \ep(1+t)\langle
	 u''(t),\roep(t)\rangle.
\end{eqnarray*}

Let us estimate the last three terms.  From (\ref{eq:tetep'}) we have
that 
$$|\roep'(t)|^{2}\leq
2\left(|\rep'(t)|^{2}+|\tetep'(t)|^{2}\right)\leq
2|\rep'(t)|^{2}+\alpha_{14}(1+t)^{-2/\ep}.$$

Moreover from (\ref{est:cep-c}) and (\ref{est:sup}) we deduce that
\begin{eqnarray*}
	(c(t)-\cep(t))\langle Au(t),\roep(t)\rangle & \leq & 
	|c(t)-\cep(t)|\cdot|A^{1/2}u(t)|\cdot|A^{1/2}\roep(t)|\\
	 & \leq & \au{\roep(t)}\cdot 
	 \left(\alpha_{5}\au{\roep(t)}\cdot\au{u(t)}\right) \\
	 & \leq & \frac{1}{4}\cep(t)\auq{\roep(t)}+
	 \frac{\alpha_{5}^{2}}{\cep(t)}\auq{\roep(t)}\auq{u(t)}\\
	 & \leq & \frac{1}{4}\cep(t)\auq{\roep(t)}+
	 \alpha_{15}\ep^{2}\auq{u(t)}.
\end{eqnarray*}

The term involving $u''(t)$ can be estimated more or less as in the
proof of (\ref{est:u''-int}). We obtain that
\begin{eqnarray*}
	|\ep\langle\roep(t), u''(t)\rangle| & \leq & \au{\roep(t)}\cdot\ep
	\left(\alpha_{1}|A^{1/2}u(t)|+\alpha_{2}(1+t)^{2}|A^{3/2}u(t)|\right)
	\\
	 & \leq & \frac{1}{4}\cep(t)\auq{\roep(t)}+
	 \alpha_{16}\ep^{2}\left(|A^{1/2}u(t)|^{2}+
	 (1+t)^{4}|A^{3/2}u(t)|^{2}
	 \right).
\end{eqnarray*}

Replacing all these estimates in the expression for $\mdep'(t)$ we
find that
\begin{eqnarray*}
	\mdep'(t) & \leq & -\frac{1}{2}(1+t)\cep(t)\auq{\roep(t)}+
	2\ep(1+t)|\rep'(t)|^{2}+\alpha_{14}\ep(1+t)^{1-2/\ep}+  \\
	 &  & +\alpha_{17}\ep^{2}(1+t)|A^{1/2}u(t)|^{2}+
	 \alpha_{16}\ep^{2}(1+t)^{5}|A^{3/2}u(t)|^{2}.
\end{eqnarray*} 

Now we integrate in $[0,t]$.  Since $\mdep(0)=0$, by
(\ref{th:parab-int}) with $k=0$ and $k=2$ we obtain that
\begin{equation}
	\frac{1}{2}\int_{0}^{t}(1+s)\cep(s)\auq{\roep(s)}\,ds\leq
	2\ep\int_{0}^{t}(1+s)|\rep'(s)|^{2}\,ds+\alpha_{18}\ep^{2}-
	\mdep(t).
	\label{est:mdep-interm}
\end{equation}

Let us estimate the last term. By the smallness of $\ep$ and
(\ref{eq:tetep'}) we have that
\begin{eqnarray*}
	-\mdep(t) & \leq & -\frac{1}{4}|\roep(t)|^{2}+
	|\roep(t)|\cdot\ep(1+t)|\roep'(t)|  \\
	 & \leq & -\frac{1}{8}|\roep(t)|^{2}+4\ep^{2}(1+t)^{2}
	 \left(|\rep'(t)|^{2}+|\tetep'(t)|^{2}\right)\\
	 & \leq & -\frac{1}{8}|\roep(t)|^{2}+
	 4\ep^{2}(1+t)^{2}|\rep'(t)|^{2}+\alpha_{19}\ep^{2}. 
\end{eqnarray*}

Going back to (\ref{est:mdep-interm}) this implies that
\begin{eqnarray}
	\lefteqn{\hspace{-2em}
	|\roep(t)|^{2}+\int_{0}^{t}(1+s)\cep(s)\auq{\roep(s)}\,ds\ \leq} 
	\nonumber \\
	 & \leq & 16\ep\int_{0}^{t}(1+s)|\rep'(s)|^{2}\,ds+
	 32\ep^{2}(1+t)^{2}|\rep'(t)|^{2}\,ds+\alpha_{20}\ep^{2}.
	\label{est:mdep-1}
\end{eqnarray}

\subparagraph{\textmd{\emph{Decay estimates on $\meep$}}}

We improve (\ref{est:meep-1}) and (\ref{est:mdep-1}) by showing that
there exists a constant $K_{2}$ such that
\begin{equation}
	(1+t)^{2}\meep(t)+\int_{0}^{+\infty}(1+s)|\rep'(s)|^{2}\,ds
	\leq K_{2}\ep^{2},
	\quad\quad\forall t\geq 0,
	\label{est:meep-final}
\end{equation}
\begin{equation}
	|\roep(t)|^{2}+\int_{0}^{t}(1+s)\auq{\roep(s)}\,ds\leq K_{2}\ep^{2}
	\quad\quad\forall t\geq 0.
	\label{est:mdep-final}
\end{equation}

These estimates imply (\ref{th:ps-u}) through (\ref{th:ps-u'int}).

From (\ref{est:meep'}) we have that
\begin{eqnarray*}
	\left[(1+t)^{2}\meep(t)\right]' & = &
	(1+t)^{2}\meep'(t)+2(1+t)\meep(t)  \\
	 & \leq & -(1-2\ep)(1+t)\frac{|\rep'(t)|^{2}}{\cep(t)}+
	 2(1+t)\auq{\roep(t)}+\\
	 &  & +(1+t)^{3}\left(\alpha_{6}|Au(t)|^{2}\auq{\roep(t)}+
	 \alpha_{7}\ep^{2}|u''(t)|^{2}\right)+ \\
	 & & + \alpha_{8}(1+t)^{2-1/\ep}\au{\roep(t)}.
\end{eqnarray*}

By the smallness assumptions on $\ep$ and estimate (\ref{est:sup})
this implies that
\begin{eqnarray*}
	\left[(1+t)^{2}\meep(t)\right]' & \leq & -\frac{1}{2}(1+t)
	\frac{|\rep'(t)|^{2}}{\cep(t)}+
	\frac{2}{\mu_{1}}(1+t)\cep(t)\auq{\roep(t)}+ \\
	 & & +\alpha_{21}\ep^{2}(1+t)^{3}\left(|Au(t)|^{2}+|u''(t)|^{2}\right)
	 +\alpha_{22}\ep(1+t)^{2-1/\ep}.
\end{eqnarray*}

Now we integrate in $[0,t]$.  By (\ref{est:mdep-1}),
(\ref{th:parab-int}) with $k=1$, and (\ref{est:u''-int}) we obtain
that
\begin{eqnarray*}
	\lefteqn{\hspace{-3em}
	(1+t)^{2}\meep(t)+\frac{1}{2}\int_{0}^{t}(1+s)
	\frac{|\rep'(s)|^{2}}{\cep(s)}\,ds\;\leq} \\
	 & \leq &
	\alpha_{23}\ep^{2}+\frac{32}{\mu_{1}}\ep\int_{0}^{t}(1+s)|\rep'(s)|^{2}\,ds+
	\frac{64}{\mu_{1}}\ep^{2}(1+t)^{2}|\rep'(t)|^{2}.
\end{eqnarray*}

Rearranging the terms this may be rewritten as
\begin{eqnarray}
	\lefteqn{\hspace{-23em}
	\left(\frac{1}{\cep(t)}-
	\frac{64}{\mu_{1}}\ep\right)(1+t)^{2}\ep|\rep'(t)|^{2}+
	(1+t)^{2}\auq{\roep(t)}+}
	\nonumber\\
	+\frac{1}{2}\int_{0}^{t}
	\left(\frac{1}{\cep(t)}-\frac{64}{\mu_{1}}\ep\right)
	(1+s)|\rep'(s)|^{2}\,ds & \leq & \alpha_{22}\ep^{2}.
	\label{est:meep-2}
\end{eqnarray}

By (\ref{defn:ep0-bis}) we have that
$$\frac{1}{\cep(t)}-\frac{64}{\mu_{1}}\ep\geq
\frac{1}{\mu_{2}}-\frac{64}{\mu_{1}}\ep\geq
\frac{1}{2\mu_{2}},$$
and therefore (\ref{est:meep-final}) easily follows from
(\ref{est:meep-2}).

This in turn allows to estimate the first two terms in the right-hand
side of (\ref{est:mdep-1}).  We thus obtain (\ref{est:mdep-final}).

\subsection{Decay-error estimates for more regular data}\label{sec:rmk}

Let us assume that $(u_{0},u_{1})\in D(A^{(k+2)/2})\times D(A^{k/2})$.
Let us consider the energies
$$\mdepk(t):=\textstyle{\frac{1-(2k-1)\ep}{2}}
(1+t)^{2k-2}|A^{(k-1)/2}\roep(t)|^{2}+
\ep(1+t)^{2k-1}\langle A^{(k-1)/2}\roep'(t),A^{(k-1)/2}\roep(t)\rangle,$$
$$\meepk(t):=\ep\frac{|A^{(k-1)/2}\rep'(t)|^{2}}{\cep(t)}+
|A^{k/2}\roep(t)|^{2},$$
which are the natural extensions of those defined in (\ref{defn:mdep})
and (\ref{defn:meep}).  For simplicity in this section $C_{k}$ denotes
a constant, which may be different from line to line, but always
depends only on $k$, on $\mu_{1}$, $\mu_{2}$, $L$, and on the norms of
initial data in the appropriate spaces (so that $C_{k}$ doesn't depend
on $\ep$ and $t$).

Working with $\mdepk(t)$ and $(1+t)^{2k}\meepk(t)$ as we did with
$\mdep(t)$ and $(1+t)^{2}\meep(t)$, with an easy induction we obtain
(\ref{th:rmk-2}) and
\begin{equation}
	(1+t)^{k}\left(\sqrt{\ep}|A^{(k-1)/2}\rep'(t)|+
	|A^{k/2}\roep(t)|\right)\leq C_{k}\ep.
	\label{th:rmk-1/2}
\end{equation}

In order to prove (\ref{th:rmk-1}) it remains to show that
\begin{equation}
	(1+t)^{k-1}|A^{(k-2)/2}\rep'(t)|\leq C_{k}\ep.
	\label{th:rmk-3}
\end{equation}

To this end we consider the energy
$$\mGepk(t):=(1+t)^{2k-2}|A^{(k-2)/2}\rep'(t)|^{2},$$
which is the natural extension of (\ref{defn:G}).  Its time derivative
is (for shortness' sake we omit the dependence on $t$)
\begin{eqnarray*}
	\mGepk' & = &
	\left(2k-2-\frac{2}{\ep}\right)(1+t)^{2k-3}|A^{(k-2)/2}\rep'|^{2}+
	\\
	 &  & -\frac{2}{\ep}(1+t)^{2k-2}\langle
	A^{(k-2)/2}\rep',\cep A^{k/2}\roep+
	(\cep-c)A^{k/2}u+
	\ep A^{(k-2)/2}u''\rangle  \\
	 & \leq &
	 \left(2k-2-\frac{2}{\ep}\right)\frac{1}{1+t}\mGepk+
	 \frac{2}{\ep}\frac{1}{1+t}\sqrt{\mGepk}\times  \\
	 & & \times(1+t)^{k}\left\{\cep|A^{k/2}\roep|+|\cep-c||A^{k/2}u|+
	 \ep|A^{(k-2)/2}u''|\right\}.
\end{eqnarray*}

All the terms in the last line can be easily estimated.  From
(\ref{th:rmk-1/2}), (\ref{th:parab-local}), (\ref{est:cep-c}), and
(\ref{est:sup}) we have indeed that 
$$(1+t)^{k}|A^{k/2}\roep(t)|\leq
C_{k}\ep,
\hspace{3em}
|\cep(t)-c(t)|(1+t)^{k}|A^{k/2}u(t)|\leq C_{k}\ep,$$
while arguing as in the proof of (\ref{est:u''-int}) we obtain that
$$(1+t)^{k}|A^{(k-2)/2}u''(t)|\leq\alpha_{1}(1+t)^{k}|A^{k/2}u(t)|+
\alpha_{2}(1+t)^{k+2}|A^{(k+2)/2}u(t)|\leq C_{k}.$$

Therefore for $\ep$ small enough it turns out that
$$\mGepk'(t)\leq\frac{1}{\ep(1+t)}\left(
-\mGepk(t)+\ep C_{k}\sqrt{\mGepk(t)}\right)
\quad\quad\forall t\geq 0.$$

Since $\mGepk(0)=0$, from Lemma~\ref{lemma:ode-sqrt} we deduce that
$\mGepk(t)\leq C_{k}^{2}\ep^{2}$ for every $t\geq 0$, which is
equivalent to (\ref{th:rmk-3}).

\subsection{The supercritical case}\label{sec:p>1}

Let us prove Theorem~\ref{thm:p>1}. Let us define $H(t)$ as in
(\ref{defn:H}). Then
$$H'(t)=-2b(t)|\uep'(t)|^{2}\geq -\frac{2}{\ep}b(t)H(t),$$
hence
$$H(t)\geq H(0)\exp\left(-\frac{2}{\ep}\int_{0}^{t}b(s)\,ds\right).$$

The right-hand side is greater than a positive constant independent on
$t$ because of (\ref{hp:int-b}) and the fact that $H(0)>0$. This
implies (\ref{th:p>1}).

\subsection*{Acknowledgements}

We would like to thank Professor Taeko Yamazaki for sending us a
preliminary version of references \cite{yamazaki-wd} and
\cite{yamazaki-cwd}, and for pointing out reference \cite{wirth}.

\label{NumeroPagine}

\end{document}